\documentclass[12pt, a4paper]{article}

\usepackage{a4wide}
\usepackage{amsmath}
\usepackage{amsfonts}
\usepackage{enumerate}
\date{}
\begin{document}
\title{\large\textbf{A NONSEPERABLE INVARIANT EXTENSION OF LEBESGUE MEASURE - A GENERALIZED AND ABSTRACT APPROACH }}
\author{\normalsize
\textbf{S.BASU \& D.SEN}\\
}
\date{}	
\maketitle
{\small \noindent \textbf{{\textbf{AMS subject classification:}}}} {\small $28A05$, $28D05$, $28A99$, $28C10$.\\}
{\small\noindent \textbf{\textbf{\textmd{\textbf{{Key words and phrases :}}}}} Transformation groups, invariant measures, nonseparable invariant extension, $Ulam(k,k^{^{+}})$-matrix, $k$-additive measurable structure, $k$-independent (strictly $k$-independent)family, almost invariant set, $k$-small system, admissible $k$-additive algebra.}\\
\vspace{.03cm}\\
{\normalsize
\textbf{\textbf{ABSTRACT:}} Here using some methods of combinatorial set theory, particularly the ones related to the construction of independent families of sets and some modified version of the notion of small sets originally introduced by Rie$\check{c}$an, Rie$\check{c}$an and Neubrunn, we give abstract and generalized formulation of a remarkable theorem of Kakutani and Oxtoby relating to nonseparable extension of Lebesgue measure in spaces with transformation groups.\\

\section{\large{INTRODUCTION}}
\normalsize\hspace{.9cm} Let $(E,\mathcal S,\mu)$ be any measure space where $E$ is a nonempty basic set, $\mathcal S$ is a $\sigma$-algebra of subsets of $E$ containing all singletons and $\mu$ a nonzero, $\sigma$-finite, diffused (or, continuous) measure on $\mathcal S$. The general measure extension problem is concerned about extending $\mu$ to a maximally large class of subsets of $E$. According to famous theorem of Ulam $[19]$, it is consistent with the Axioms of set theory that so long as no cardinal less than card\hspace{.05cm}{$E$} is an weakly inaccessible cardinal, there can be no extension of $\mu$ which coincides with the power set of $E$. Thus there always exists $X\subseteq E$ such that $X\notin dom (\mu)$ and $\mu$ can be further extended to $\mu^{\prime}$ with $X\in dom (\mu^{\prime})$. Therefore, under the assumption that there is no large cardinal, we may infer that no maximal extension of $\mu$ is possible.\\
\vspace{.1cm}

\hspace{.3cm}  Investigations into the extension problems in measure theory have useful applications in many other branches of modern mathematics such as axiomatic set theory, general topology, functional analysis, probability theory, etc and possibly started with the work of Szpilrajn (Marczewski). In his classical papers $[17]$, $[18]$, he gave several such constructions of invariant extensions and also added to it a list of problems. His method is sufficiently general so as to make it applicable for any complete measure space and is based upon the idea of adding to the $\sigma$-ideal $\mathcal I(\mu)$ of $\mu$-null sets, some new sets which are $\mu$-nonmeasurable and having inner $\mu$-measure zero. This yields the $\sigma$-ideal $\mathcal I^{^{\prime}}$ and in turn generates the $\sigma$-algebra $\mathcal S^{^{\prime}}$ from $\mathcal S$ and $\mathcal I^{^{\prime}}$. Any arbitrarily chosen element $U$ of $\mathcal S^{^{\prime}}$ admits a representation
of the form $U = (X\setminus Y)\cup Z$ where $X\in \mathcal S$, and $Y,Z\in \mathcal I^{^{\prime}}$ and on $\mathcal S^{^{\prime}}$ we define a measure $\mu^{\prime}$ by putting $\mu^{\prime}(U) = \mu(X)$. It may be observed that this way of extending $\mu$ do not change its metrical structure and the two metric spaces associated with $\mu$ and $\mu^{\prime}$ have the same topological weights. Therefore, if the original measure $\mu$ is separable (i.e the metric associated with $\mu$ is separable) then so is $\mu^{\prime}$.\\
\vspace{.02cm}

Many years ago, it was asked whether there is some nonseparable invariant extension of Lebesgue measure in $\mathbb{\textbf{R}}$ and researches in this directions were carried out by Kakutani and Oxtoby $[10]$, Kodaira and Kakutani $[9]$ and several others. The extension obtained by Kakutani and Oxtoby has character $2^{^{c}}$ and that by Kodaira and Kakutani has character $c$ and both are invariant under the groups of all isometries in $\mathbb{\textbf{R}}$. Moreover, the method of Kakutani and Oxtoby can be generalized for Haar measure in infinite compact metrizable groups $[2]$.\\
\vspace{.02cm}

Let $(\Omega,\mathcal L, \lambda)$ be the Lebesgue measure space where $\Omega = [0,1)$, $\mathcal L$ is the class of Lebesgue measurable sets in $\Omega$, and $\lambda$ the restriction of the usual Lebesgue measure on $\Omega$. We call a measure space $(\Omega, \mathcal M, m)$ an extension of $(\Omega, \mathcal L, \lambda)$ $[10]$ if $\mathcal L \subseteq \mathcal M$ and $m(M) = \lambda(M)$ for every $M$$\in\mathcal L$, and, it is a proper extension if $\mathcal L\neq\mathcal M$. The character $[10]$ of any extension $(\Omega, \mathcal M, m)$ is the smallest cardinal for which there exists a subfamily $\mathcal U$ of $\mathcal M$ of that cardinality and having the property that for every $M\in \mathcal M$ and every $\epsilon > 0$, there exists $A\in \mathcal U$ such that $m(M \Delta A)< \epsilon$.\\
\vspace{.02cm}

The extension is called separable $[10]$ if card\hspace{.1cm}($\mathcal U$)\hspace{.03cm}$\leq\omega_{_{0}}$ where $\omega_{_{0}}$ is the ordinal corresponding to the first infinite cardinal. An one-to-one mapping T : $\Omega\rightarrow \Omega$ is called $\mathcal M-m$ preserving $[10]$ if $M\in\mathcal M$ implies that T$(M)$$\in\mathcal M$, T$^{^{-1}}$$(M)\in\mathcal M$ and $m($T$(M))$ = $m($T$^{^{-1}}(M)$ = $m(M)$.\\
\vspace{.1cm}

\hspace{.3cm} In their remarkable paper $[10]$, Kakutani and Oxtoby proved the following theorem which serves as a motivation for the present work. Here, we generalize this theorem through some abstract approach, where, instead of measure, more delicate structures in spaces with transformation groups are used.\\
\vspace{.2cm}\\
\textbf{THEOREM(KO):}\textit{ There exists an extension $(\Omega, \mathcal M, m)$ of the Lebesgue measure space $(\Omega, \mathcal L, \lambda)$ having the following properties:\\
(i) The character of $(\Omega, \mathcal M, m)$ is $2^{^{c}}$.\\
(ii) Every $\mathcal L-\lambda$ preserving transformation is also $\mathcal M-m$ preserving.}\
\vspace{.001cm}
\section{\large{PRELIMINARIES AND RESULTS}}
\normalsize\hspace{.9cm} For our purpose, we use some methods of combinatorial set theory, particularly the ones that are related to the construction of independent families of sets in infinite base space and also some modified version of the notion of small sets. It is worthwhile to note here that the notion of small sets (or, families) was originally introduced by Rie$\check{c}$an, Rie$\check{c}$an and Neubrunn $[13]$, $[14]$ (see also $[15]$, $[16]$)while giving abstract formulations of some well known classical theorems on measure and integration.\\
\vspace{.1cm}

\hspace{.3cm} By a space $X$ equipped with a transformation group $G$ we mean a pair $(X,G)$ where $X$ is a nonempty set and $G$ is a subgroup of the symmetric group Symm$(X)$ of all bijections from $X$ onto $X$ satisfying the following two conditions:\\
\textit{(i)} for each $g\in G$, $x\rightarrow gx$ is a bijection (or, permutation) of $X$.\\
\textit{(ii)} for all $x\in X$ and $g_{_{1}}, g_{_{2}} \in G$, $g_{_{1}}(g_{_{2}}x) = g_{_{1}}g_{_{2}}x$.\\
\vspace{.1cm}

\hspace{.2cm} We say that $G$ acts freely on $X$ if $\{x \in X : gx = x\} = \emptyset$ for all $g\in G\setminus \{e\}$ where $`e'$ is the identity element of $G$ (in fact, $e : X\rightarrow X$ is the identity transformation on X)and $G$ acts transitively on $X$ if $Gx=X$ for every $x\in X$. For any $g \in G$ and $E \subseteq X$, we write $gE$ to denote the set $\{gx : x \in E\}$ and call a nonempty family (or, class) $\mathcal A$ of sets as $G$-invariant $[3]$ if $gE\in\mathcal A$ for every $g\in G$ and $E\in\mathcal A$. If $\mathcal A$ is a $\sigma$-algebra, then a measure $\mu$ on $\mathcal A$ is called $G$-invariant $[3]$ if $\mathcal A$ is a $G$-invariant class and $\mu(gE) = \mu(E)$ for every $g\in G$ and $E\in \mathcal A$. It is called $G$-quasiinvariant $[3]$ if $\mathcal A$ and the $\sigma$-ideal generated by $\mu$-null sets are both $G$-invariant classes. Obviously, any $G$-invariant measure is also $G$-quasiinvariant but not conversely. Any set of the form $Gx = \{gx : g\in G\}$ for some $x\in X$ is called a $G$-orbit of $x$. The collection of all such $G$-orbits give rise to a partition of $X$ into mutually disjoint sets. A subset $E$ of $X$ is called a complete $G$-selector (or, simply, a $G$-selector) in $X$ if $E\hspace{.0005cm}\cap\hspace{.0005cm}Gx$ consists of exactly one point for each $x\in X$.\

Throughout this paper, we identify every infinite cardinal with the least ordinal representing it, and, every ordinal (in the sense introduced by Von Neumann) with the set of all ordinals preceeding it. We write card $A$ and card $\mathcal A$ to denote the cardinals of any set $A$ or any class $\mathcal A$ of sets and as is usually done elsewhere, express the first infinite and the first uncountable cardinals by the notations $\omega_{_{0}}$ and $\omega_{_{1}}$ respectively. For any infinite cardinal, we use symbols like $\xi$, $\varrho$, $\eta$, $k$ etc and write $k^{^{+}}$ for the successor of $k$. Throughout the entire discourse, we stipulate to work within the frammework of ZFC.\\
\vspace{.1cm}

Let $k$ be any infinite cardinal and suppose card $X \geq k$, we define \\
\vspace{.2cm}\\
\textbf{DEFINITION \textbf{$\textbf{2.1[1]}$:}} A\textit{ pair $( \mathcal S, \mathfrak I)$ as a $k$-additive measurable structure on $(X , G)$ if\\
(i) $\mathcal S$ is an algebra and $\mathfrak I (\subseteq \mathcal S)$ a proper ideal on $X$.\\
(ii) Both $\mathcal S$ and $\mathfrak I (\subseteq \mathcal S)$ are $k$-additive in the sense that they are closed with respect to the union of at most $k$ number of sets, and\\
(iii) $\mathcal S$ and $\mathfrak I$ are $G$-invariant classess}.\

Henceforth, by a $k$-additive algebra (resp. ideal) on $(X, G)$ we will mean that it is a $k$-additive algebra (resp. ideal) on $X$ and also $G$-invariant. In particular, if $G$ consists only of the identity transformation on $X$, then $(\mathcal S, \mathfrak I)$ is simply called a $k$-additive measurable structure on $X$.\\
\vspace{.2cm}\\
\textbf{DEFINITION \textbf{$\textbf{2.2[1]}$:}} A \textit{ $k$-additive measurable structure $(\mathcal S, \mathfrak I)$ on $(X , G)$ is called $k^{^{+}}$-saturated if the cardinality of any arbitrary collection of mutually disjoint sets from $\mathcal S\setminus \mathfrak I$ is atmost $k$}.\

The notion of $\omega_{_{0}}$-additive measurable structure on a nonempty basic set $E$ was defined by Kharazishvili. In $[4]$, this was referred to as a measurable structure consisting of a pair $(\mathcal S, \mathcal I)$ where $\mathcal S$ is a $\sigma$-algebra and ${\mathcal I \subseteq \mathcal S}$ a proper $\sigma$-ideal of sets in $E$. If $E$ is a group and $\mathcal S$, $\mathcal I$ are $G$-invariant classes, then $(\mathcal S, \mathcal I)$ according to Kharazishvili is a $G$-invariant measurable structure on $E$. Using the notion of a measurable structure as introduced by him, Kharazishvili proved several interesting results in commutative (or, more generally in solvable) groups $[4]$. In $[5]$, he used similar type of structures to generalize two classical results of Sierspinski. It may be noted that the notion of a $k$-additive, $k^{^{+}}$-saturated measurable structure on $(X , G)$ lies in between a $k$-additive measurable structure on $(X , G)$ satisfying countable chain condition (or, Suslin condition) and a $\omega_{_{0}}$-additive measurable structure which is $k^{^{+}}$-saturated.\\
\vspace{.2cm}\\
\textbf{DEFINITION \textbf{$\textbf{2.3[1]}$:}} \textit{In a space $(X,G)$ with a transformation group $G$, a set $E\subseteq X$ is called almost $G$-invariant with respect to some ideal $\mathfrak I$ if $gE\hspace{.04cm}\Delta\hspace{.04cm}E \in \mathfrak I$ for every $g\in G$}.\

If the ideal $\mathfrak I$ is $k$-additive, then it can be easily checked that the class of all sets which are almost $G$-invariant with respect to $\mathfrak I$ forms a $k$-additive algebra on $X$. Further,\\
\vspace{.2cm}\\
\textbf{DEFINITION \textbf{$\textbf{2.4[1]}$:}} A \textit{set $E\subseteq X$ is called $(\mathcal S,\mathfrak I)$-thick if $B\subseteq {X\setminus E}$ and $B\in \mathcal S$ implies that $B\in \mathfrak I$.}\\
\vspace{.2cm}\\
\textbf{\textbf{PROPOSITION $\textbf{2.5}$ :}} \textit{Assume that the pair $(\mathcal S ,\mathfrak I)$ is a $k$-additive measurable structure on $(\textsl {X} , \textsl {G})$ which is also $k {^{^{+}}}$-saturated and let $\textsl{E} \subseteq \textsl {X}$ be a set which is almost $\textsl {G}$-invariant with respect to $\mathfrak I$. Then $\textsl{E} \in \mathcal S$ implies either $\textsl{E} \in \mathfrak I$ or $\textsl {X} \setminus \textsl{E} \in \mathfrak I$. If $ E \notin \mathfrak I$, then $\textsl{E}$ is $(\mathcal S , \mathfrak I)$-thick in $\textsl {X}$.}\\
\vspace{.1cm}\\
\textbf{Proof :} Let $\textsl{E}\in\mathcal S$. If $\textsl{E} \in \mathfrak I$, then there is nothing to prove. Suppose $\textsl{E} \notin \mathfrak I$. Then $\textsl {X} \setminus \textsl{E} \in \mathfrak I$, for otherwise, it is possible to generate by transfinite recursion a set $\{ g_{\alpha} : \alpha < k\}$ of points in $\textsl {G}$ such that $\textsl {X}\setminus{\displaystyle\bigcup_{_{_{ _{0\leq \alpha < k}}}}}{g_{\alpha}\textsl{E}} \in \mathfrak I$. But this contradicts the hypothesis.\

Now let $\textsl{E}\notin \mathcal S$. Then obviously $\textsl{E}\notin \mathfrak I$ and also $\textsl {X} \setminus \textsl{E} \notin \mathfrak I$. If $\textsl{E}$ is not $(\mathcal S , \mathfrak I)$-thick, then there should exist $\textsl{B} \in \mathcal S \setminus \mathfrak I$ such that $\textsl{B} \subseteq \textsl {X} \setminus \textsl{E}$. By a similar reasoning as given above, there exists a transfinite $k$-sequence $\{ h_{\alpha} : \alpha < k\}$ in $\textsl {G}$ such that $\textsl {X}\setminus{\displaystyle\bigcup_{_{_{ _{0\leq \alpha < k}}}}}{h_{\alpha}\textsl{B}} \in \mathfrak I$. But then by virtue of $k$-additivity of $\mathfrak I$, there exists some $\alpha _{_{0}}$$< k$ such that $E\hspace{.04cm} \cap\hspace{.04cm} {h _{\alpha}{_{_{_{0}}}} \textsl {B}}\notin \mathfrak I$ which again contradicts the hypothesis.\\
\vspace{.2cm}

\hspace{.2cm}The notion of a small system was introduced by Rie$\check{c}$an $[13]$, Rie$\check{c}$an and Neubrunn $[15]$ (see also $[14]$,$[16]$). As an initial step towards giving an abstract generalization the above theorem(Theorem(KO))of Kakutani and Oxtoby, we build up a $k$-additive measurable structure on $(X,G)$ from a $k$-additive algebra $\mathcal S$ on $(X,G)$ and a transfinite $k$-sequence which in the present situation is a modified and generalized version of systems of small sets or small systems originally introduced by Rie$\check{c}$an and Neubrunn.\\
\vspace{.2cm}

We call it a $k$-small system on $(X,G)$ and define it as\\
\vspace{.1cm}

\textbf{DEFINITION \textbf{$\textbf{2.6}$:}}
  A\textit{ transfinite $k$-sequence $\{\mathcal N_{_{\alpha}}\}_{_{0\leq \alpha< k}}$ where each $\mathcal N_{_{\alpha}}$ is a class of sets in $X$ satisfying the following set of conditions :\\
(i) $\emptyset\in \mathcal N_{_{\alpha}}$ for $\alpha < k$.\\
(ii) Each $\mathcal N_{_{\alpha}}$ is a $G$-invariant class.\\
(iii) If $E\in \mathcal N_{_{\alpha}}$ and $F\subseteq E$, then $F\in \mathcal N_{_{\alpha}}$.\\
(iv) $E\in \mathcal N_{\alpha}$ and $F\in{\displaystyle{\bigcap_{0\leq\alpha<k}}\hspace{.01cm}{\mathcal N_{\alpha}}}$ implies $E\cup F\in\mathcal N_{\alpha}$.\\
(v) For any $\alpha < k$, there exists $\alpha^{\ast}> \alpha$ such that for any one-to-one correspondence $\beta\rightarrow \mathcal N_{_{\beta}}$ with $\beta > \alpha^{\ast}$, $\cup{E_{_{\beta}}}\in \mathcal N_{_{\alpha}}$ whenever $E_{_{\beta}}\in \mathcal N_{_{\beta}}$. \\
(vi) For any $\alpha, \beta <k$, there exists $\gamma> \alpha, \beta$ such that $\mathcal N_{_{\gamma}}\subseteq\mathcal N_{_{\alpha}}$ and $\mathcal N_{_{\gamma}}\subseteq \mathcal N_{_{\beta}}$.\\
}\vspace{.001cm}

We further say that\\
\vspace{.1cm}

\textbf{DEFINITION \textbf{$\textbf{2.7}$:}} A \textit{$k$-additive algebra $\mathcal S$ on $(X,G)$ is admissible with respect to the $k$-small system $\{\mathcal N_{\alpha}\}{_{_{0\leq\alpha<k}}}$ if\\
(i) $\mathcal S\setminus\mathcal N_{\alpha}\neq\emptyset$ for some $\alpha<k$ and $\mathcal S\cap\mathcal N_{\alpha}\neq\emptyset$ for every $\alpha<k$.\\
(ii) $\mathcal N_{\alpha}$ has a $S$-base i.e every $E\in\mathcal N_{\alpha}$ is contained in some $F\in \mathcal N_{\alpha}\cap\mathcal S$,\\
and (iii) $\mathcal S\setminus\mathcal N_{\alpha}$ satisfies $k$-chain condition, i.e, the cardinality of any arbitrary collection of mutually disjoint sets from $\mathcal S\setminus\mathcal N_{\alpha}$ is atmost $k$.}\
\vspace{.01cm}

By \textit{(i)} in the above definition, we mean that $\mathcal S$ is compatible with $\{\mathcal N_{\alpha}\}{_{_{0\leq\alpha<k}}}$; by \textit{(ii)} we mean that $\mathcal S$ constitutes a base for $\{\mathcal N_{\alpha}\}{_{_{0\leq\alpha<k}}}$ and by \textit{(iii)} we express that $\mathcal S$ satisfies the $k$-chain condition with respect to $\{\mathcal N_{\alpha}\}{_{_{0\leq\alpha<k}}}$.\
\vspace{.1cm}

We set $\mathcal N_{_{\infty}} = {\displaystyle\bigcap_{_{0\leq\alpha <k}}}{\mathcal N_{_{\alpha}}}$. By virtue of conditions \textit{(ii)}, \textit{(iii)} and \textit{(v)} of Definition $2.6$, it is easy to check that $\mathcal N_{_{\infty}}$ is a $k$-additive ideal on $(X,G)$. Let $\widetilde{\mathcal S}$ be the $k$-additive algebra on $(X,G)$ generated by $\mathcal S$ and $\mathcal N_{_{\infty}}$, where each element of $\widetilde{\mathcal S}$ admits a representation of the form $(X\setminus Y)\cup Z$ where $X\in \mathcal S$ and $Y,Z\in \mathcal N_{_{\infty}}$. Thus$(\widetilde{\mathcal S}, \mathcal N_{_{\infty}})$ is a $k$-additive measurable structure on $(X,G)$. Moreover,\\
\vspace{.2cm}

\textbf{THEOREM $\textbf{2.8:}$} \textit{If $\mathcal S$ is admissible with respect to $\{\mathcal N_{\alpha}\}{_{_{0\leq\alpha<k}}}$, then $(\mathcal{\widetilde S}, \mathcal N_{\infty})$ is $k^{^{+}}$-saturated.}\\
\vspace{.2cm}

 The existence of an independent family of subsets of an infinite set, with maximal cardinality was solved by Tarski $[11]$. He showed that such a family exists and has cardinality ${{2}}^{^{^{card \textsl{(E)}}}}$. The result has many interesting applications. One such is its use in proving that the cardinality of all ultrafilters defined on an arbitrary infinite set $E$ is $2 {^{^{ {2}}}}^{^{^{card \textsl {(E)}}}}$. However, if the cardinality of the basic set is that of the continuum $c$, then the existence of a strictly independent family of subsets of $E$ having cardinality $2 {^{^{c}}}$ can be proved and this result has an application in the construction of a nonseparable invariant extension of the Lebesgue measure space $[2]$. Below, we introduce a more general definition for a $k$-independent (resp. strictly $k$-independent) family where $k$ is an arbitrary infinite cardinal.\\
\vspace{1cm}\\
%\newpage
\textbf{Definition $\textbf{2.9}$ :} A\textit{ family $\{ \textsl{A} _{i} : i \in \textsc{I}\}$ of subsets of $\textsl {X}$ is called $k$-independent (resp. strictly $k$-independent) if for each set ${\textsc{J}\subseteq \textsc{I}}$ having card$\hspace{.07cm}$\textsc{J}${<k}$ (resp. card$\hspace{.07cm}$\textsc{J}${\leq k}$) and every function $f : \textsc{J} \rightarrow \{0,1\}$, we have $\cap \{ \textsl{A}_{j}^{f(j)} : j \in \textsc{J} \}\neq \emptyset$ where $\textsl{A}_{j}^{f(j)} = \textsl{A} _{j}$ if $f(j) = 0$ and $\textsl{A}_{j}^{f(j)} = \textsl {X} \setminus \textsl{A} _{j}$ if $f(j) = 1$.}\\
\vspace{.1cm}

In particular, the definition of an independent or $\omega_{_{0}}$-independent (in the set theoretic sense) family is already given in $[6]$. The above general definition is frammed in this pattern. For another introduction to $k$-independent (resp. strictly $k$-independent) family, see $[12]$. However, Definition $2.9$ can be further generalized using Definition $2.1$.\\
\vspace{.1cm}\\
\textbf{Definition $\textbf{2.10}$ :} A\textit{ family $\{\textsl{A} _{i} : i \in \textsc{I} \}$ of subsets of $\textsl {X}$ is $k$-independent (resp. strictly $k$-independent) with respect to any $k$-additive measurable structure $(\mathcal S , \mathfrak I)$ on $(\textsl {X} , \textsl {G})$ if for each set ${\textsc{J}\subseteq\textsc{I}}$ having card$\hspace{.07cm}$\textsc{J}${<k}$ (resp. card$\hspace{.07cm}$\textsc{J}${\leq k}$) and every function $ f : \textsc{J} \rightarrow \{0 , 1\}$, $ \textsl{B }\subseteq \textsl {X} \setminus \cap \{ \textsl{A} _{j} ^{f(j)} : j \in \textsc{J} \}$ and $ \textsl{B}\in\mathcal S$ implies that $ \textsl{B} \in \mathfrak I$, where $\textsl{A} _{j} ^{f(j)}$ ($ j \in \textsc{J}) $ has the same meaning as before.}\

Note that in the above Definition, condition\textit{(iii)} of Definition $2.1$ plays no role. So we may think of it as a $k$-independent(resp. strictly $k$-independent) family with respect to some $k$-additive measurable structure $(\mathcal S,\mathfrak I)$ on $X$. The notion of an independent (resp. strictly independent) family with respect to a measure was earlier given in $[7]$.So the above Definition is just an extension of this concept given in terms of some $k$-additive measurable structure and viewed further in the context of $(\mathcal S,\mathfrak I)$-thick sets (Definition $2.4$), we may call a family $\{A_{i} : i \in\textsc{I}\}$ $k$-independent (resp. strictly $k$-independent) with respect to $(\mathcal S,\mathfrak I)$ on $X$ if for each ${\textsc{J}\subseteq\textsc{I}}$ having card$\hspace{.07cm}$\textsc{J}${<k}$ (resp. card$\hspace{.07cm}$\textsc{J}${\leq k}$) and each function $ f : \textsc{J} \rightarrow \{0 , 1\}$, the set $\cap \{A _{j}^{f(j)} : j \in\textsc{J}\}$ is $(\mathcal S,\mathfrak I)$-thick in $X$.\\
\vspace{.1cm}

We now establish, under the assumption of generalized continuum hypothesis, the existence of a family of sets in $X$ which is strictly $k$-independent with respect to the $k$-additive measurable structure $(\widetilde{\mathcal S}, \mathcal N_{_{\infty}})$ on $(X,G)$.\\
\vspace{.2cm}\\
\textbf{\textbf{PROPOSITION $\textbf{2.11}$[12] :}} \textit{Assume that the generalized continuum hypothesis holds. Then for any two infinite cardinals $\lambda$, $k$ where $\lambda$$<k$, we have $k^{^{\lambda}} = k$  provided $\lambda$ is not cofinal with $k$.}\\
\vspace{.18cm}\\
\textbf{\textbf{PROPOSITION $\textbf{2.12}$[12] :}}\textit{ Let $\textsl{E}$ be an infinite set satisfying the condition card $\textsl{E} ^{\hspace{.02cm}}{^{k}} = card\hspace{.1cm} \textsl{E}$, where $k$ is an infinite cardinal. Then there exists a maximal strictly $k$- independent family $\{\textsl{A} _{i} : i \in \textsc{I} \}$ of subsets of $\textsl{E}$ such that card $\textsc{(I)} = 2 {^{^{card \textsl{(E)}}}}$.}\\
\vspace{.1cm}\\
\textbf{\textbf{THEOREM $\textbf{2.13}$ :}} \textit{Let $(X,G)$ be a space with transformation group $G$ where card($G$)=$k^{^{+}}\leq$ card($X$). Suppose $G$ acts freely on $X$ and there exists a $G$-selector $L\in \mathcal S$. Further assume that there exist a $k$-additive algebra $\mathcal S$ and a $k$-small system $\{\mathcal N_{\alpha}\}{_{_{\alpha<k}}}$ on ($X,G$) such that $\mathcal S$ is admissible with respect to $\{\mathcal N_{\alpha}\}{_{_{\alpha<k}}}$. Then under the assumption of generalized continuum hypothesis, there exists a family $\{\textsl{A} _{i} : i \in \textsc{I} \}$ of sets in $X$ with card $\textsc{(I)} = 2^{k^{^{+}}}$ which is strictly $k$-independent with respect to $(\widetilde{\mathcal S}, \mathcal N_{_{\infty}})$}.\\
\vspace{.1cm}\\
\textbf{Proof :} We write $G$ in the form $G = {\displaystyle\cup_{_{\xi < k^{^{+}}}}} {G_{_{\xi}}}$, where $\{G_{_{\xi}} : \xi < k^{^{+}}\}$ is an increasing family of subgroups of $G$ satisfying $G_{_{\xi}} \neq \cup_{_{\eta < \xi}}{G_{_{\eta}}}$ and card $(G_{_{\xi}})$$\leq k$ for every $\xi$$< k^{^{+}}$ (for the above representation, see [6], Exercise $19$, Ch $3$).\\
\vspace{.001cm}\

\hspace{.5cm} Since $G$ acts freely on $X$, the above increasing family yields a disjoint covering $\{\Omega_{_{\gamma}} : \gamma < k^{^{+}}\}$ of $X$ where $\Omega_{_{\gamma}} = (G_{_{\gamma}} \setminus {\cup_{_{\eta < \gamma}}{G_{_{\eta}}}})L$. Moreover, as $L\in \mathcal S$, $G$ acts freely on $X$ and $\mathcal S$ is admissible with respect to $\{\mathcal N_{\alpha}\}{_{_{\alpha<k}}}$, so  $(\widetilde{\mathcal S}, \mathcal N_{_{\infty}})$} is $k^{^{+}}$-saturated (Theorem $2.8$). Therefore $gL\in \mathcal N_{_{\infty}}$ for every $g\in G$.\\
\vspace{.001cm}\

\hspace{.1cm} Now consider the $Ulam (k,k^{^{+}})$-matrix $[1]$ $(\Pi_{_{\xi,\rho}})_{\xi < k, \rho < k^{^{+}}}$ over $k^{^{+}}$ and set $E_{_{\xi, \rho}} = {\cup_{_{\gamma\in \Pi_{_{\xi, \rho}}}}}{\Omega_{_{\gamma}}}$. Then there exists $\xi_{_{0}}$ and a subset $\Xi$ of $k^{^{+}}$ having card$(\hspace{.06cm}\Xi)={k^{^{+}}}$ such that $E_{_{\xi_{_{0}}}, \rho}\notin \mathcal N_{_{\infty}}$ for $\rho\in \Xi$ and are mutually disjoint. This is so because $\mathcal N_{_{\infty}}$ is $k$-additive and $X\notin \mathcal N_{_{\infty}}$. Consequently, there exists  ${\widetilde{\Xi}\subseteq \Xi}$ with card$(\hspace{.06cm}\widetilde{\Xi})=k^{^{+}}$ and $\alpha_{_{0}}$$< k$ such that $E_{_{\xi_{_{0}}}, \rho}\notin \mathcal N_{_{\alpha_{_{0}}}}$ for every $\rho\in \widetilde{\Xi}$. Moreover, each $E_{_{\xi_{_{0}}}, \rho}$ for $\rho\in \Xi$ and more generally any union of such sets is almost $G$-invariant with respect to $\mathcal N_{_{\infty}}$ which follows from the constructions of the sets $\Omega_{_{\gamma}}$.\\
\vspace{.00001cm}\

\hspace{.1cm} Now note that $k$ is not cofinal with $k^{^{+}}$. This is so because $k$ is not cofinal with $2^{^{k}}$ and $2^{^{k}}$= $k^{^{+}}$ under the assumption of generalized continuum hypothesis. Hence according to Proposition $2.11$ and Proposition $2.12$, it follows that there exists a strictly $k$-independent family $\{\widetilde{\Xi_{_{i}}} : i\in\textsc{I}\}$ of subsets of $\widetilde{\Xi}$ such that card$\hspace{.06cm}\textsc{(I)}=2^{k^{^{+}}}$. This means that for every set ${\textsc{J}\subseteq\textsc{I}}$ having card$\hspace{.07cm}$\textsc{(J)}${\leq k}$ and every function $f :\textsc{J}\rightarrow \{0,1\}$, $\cap_{_{j\in\textsc{J}}}\hspace{.01cm}{\widetilde{\Xi}_{j}^{f(j)}} \neq \emptyset$. Consequently, $\cap_{_{j\in \textsc{ J}}}\hspace{.01cm}{A_{j}^{f(j)}} \neq \emptyset$ where $A_{i} = \cup_{_{\rho \in \widetilde {\Xi}_{i}}}\hspace{.01cm}{E_{_{\xi_{_{0}},\rho}}}$ for $i\in\textsc{I}$ making $\{A_{j} : i\in\textsc{I}\}$ a strictly $k$-independent family of sets in $X$. Moreover, this family is also strictly $k$-independent and hence $k$-independent with respect to the $k$-additive measurable structure $(\widetilde{\mathcal S}, \mathcal N_{_{\infty}})$ on $(X,G)$ which follows from Proposition $2.5$ according to which each $E_{_{\xi_{_{0}}}, \rho}$ ($\rho\in\Xi$) is $(\widetilde{\mathcal S}, \mathcal N_{_{\infty}})$- thick in $X$.\

This proves the theorem.\\
\vspace{.00001cm}\

Again from Proposition $2.5$ it follows that none of the sets $A_{i}$ $(i\in\textsc{I})$ are $(\widetilde{\mathcal S}, \mathcal N_{_{\infty}})$-measurable; in otherwords, $A_{i}\notin \widetilde{\mathcal S}$ for any $i\in\textsc {I}$. Let us consider the algebra generated by $\widetilde{\mathcal S}$ and $\{A_{i} : i\in\textsc{I}\}$. Its individual members are sets of the form ${\displaystyle\bigcup_{_{f\in \{0,1\}^{^{p}}}}(E_{_{f}}\cap (\displaystyle\cap_{j=1}^{p}\hspace{.01cm}A_{j}^{f(j)}))}$, where $E_{_{f}}\in\widetilde{\mathcal S}$, $p$ is any natural number, $\{0,1\}^{^{p}}$ is the set of all functions $f : \{1,2,\dots,p\}\rightarrow \{0,1\}$ and $\textsl{A}_{j}^{f(j)} = \textsl{A} _{j}$ if $f(j) = 0$ and $\textsl{A}_{j}^{f(j)} = \textsl {X} \setminus \textsl{A} _{j}$ if $f(j) = 1$. It is not hard to verify that the above algebra and $G$-invariant. Let $\widetilde{\widetilde{\mathcal S}}$ be the $k$-additive algebra generated by $\widetilde{\mathcal S}$ and $\{A_{i} : i\in\textsc{I}\}$. We will prove next that the extension $(\widetilde{\widetilde{\mathcal S}}, \mathcal N_{_{\infty}})$ of $(\widetilde{\mathcal S}, \mathcal N_{_{\infty}})$ serves to yield the desired generalization of Theorem(KO).\\
\vspace{.00001cm}\

We already note from the introduction that a nonseparable extension is that the cardinality of whose character is greater than $\omega_{_{0}}$. In any measure space $(X,\mathcal S, \mu)$, there is an alternative way to formulate this phenomenon. It is done by using the metrical structure of $\mu$. It is an well known fact $[8]$ that the following two assertions are equivalent:\\
(a) $\mu$ is nonseparable\\
(b) there exists an $\epsilon > 0$ and a family $\mathcal H(\epsilon) = \{Y : Y\in \mathcal S\}$ having card $\mathcal H(\epsilon) > \omega_{_{0}}$ such that for the corresponding metric $d$, $d(Y,Z) \geq \epsilon$ for $Y,Z\in \mathcal H(\epsilon)$ and $Y\neq Z$. where $d(Y,Z)=\mu(Y\Delta Z)$.\

Since there is no measure space structure available to us, we circumvent this problem by developing an uniform structure on $\widetilde{\widetilde{\mathcal S}}$. This structure is a more modified form of an uniform structure already used in $[14]$.\\
\vspace{.1cm}\\

 Let $\mathcal V = \{ U_{_{\alpha}}\}_{_{\alpha <k}}$ be a $k$-sequence defined by setting $U_{_{\alpha}} = \{(E,F)\in \widetilde{\widetilde{\mathcal S}}\times \widetilde{\widetilde{\mathcal S}} : E\Delta F\in \mathcal N_{_{\alpha}}\}$. Since $\emptyset \in \mathcal N_{_{\alpha}}$ for $\alpha < k$, so the diagonal is contained in every member of $\mathcal V$. Also, it is evident from the definition of $\mathcal V$ that $U_{_{\alpha}}^{^{-1}}\in \mathcal V$ whenever $U_{_{\alpha}}\in \mathcal V$.\\
\vspace{.00001cm}\

 From conditions \textit{(v)} and \textit{(vi)} of Definition $2.6$ for any $\beta, \gamma > \alpha^{^{\ast}}$, $\mathcal N_{_{\beta}} \cup \mathcal N_{_{\gamma}}\subseteq \mathcal N_{_{\alpha}}$ and there exists $\delta$ such that $\mathcal N_{_{\delta}}\subseteq \mathcal N_{_{\beta}}, \mathcal N_{_{\gamma}}$. Hence, $U_{_{\delta}}\ast U_{_{\delta}} = \{(E,F)\in \widetilde{\widetilde{\mathcal S}}\times \widetilde{\widetilde{\mathcal S}} : there\hspace{.01cm}\ exists\ G\in \widetilde{\widetilde{\mathcal S}}\ such \ that\ (E,G)\in U_{_{\delta}}, (G,F)\in U_{_{\delta}}\}$\\
\hspace{.1cm} $= \{(E,F)\in \widetilde{\widetilde{\mathcal S}}\times \widetilde{\widetilde{\mathcal S}} : there\hspace{.01cm}\ exists\ G\in \widetilde{\widetilde{\mathcal S}}\ such \ that\ E \Delta G\in \mathcal N_{_{\delta}}, G\Delta F\in \mathcal N_{_{\delta}}\}$\\
$\subseteq \{(E,F)\in \widetilde{\widetilde{\mathcal S}}\times \widetilde{\widetilde{\mathcal S}} : there\hspace{.01cm}\ exists\ G\in \widetilde{\widetilde{\mathcal S}}\ such \ that\ E\Delta F \subseteq (E\Delta G)\Delta (G\Delta F) \subseteq (E\Delta G)\cup (G\Delta F)\in \mathcal N_{_{\beta}}\cup \mathcal N_{_{\gamma}} \subseteq \mathcal N_{_{\alpha}}\}\subseteq U_{_{\alpha}}$. Also, from condition \textit{(vi)} of Definition $2.6$ it follows that for every $U_{_{\alpha}}$, $U_{_{\beta}}$ there exists $U_{_{\gamma}}$ such that $U_{_{\gamma}}\subseteq U_{_{\alpha}}\cap U_{_{\beta}}$.\

Hence $\mathcal V = \{ U_{_{\alpha}}\}_{_{\alpha <k}}$ forms the base of some uniformity on $\widetilde{\widetilde{\mathcal S}}$.\
\vspace{.2cm}

Thus we draw the conclusion \\
\vspace{.0001cm}\\
\textbf{THEOREM }$\textbf{2.14}$ : \textit{Let $(X,G)$ be a space with transformation group where card($X$)=$k^{^{+}}$ and $G$ acts freely and transitively on $X$. Let $\mathcal S$ be a diffused, $k$-additive algebra on $(X,G)$ which is admissible with respect to a $k$-small system $\{\mathcal N_{\alpha}\}_{_{\alpha<k}}$. Then under the assumption of generalised continuum hypothesis, we can extend the $k$-additive measurable structure $(\widetilde{\mathcal S}, \mathcal N_{_{\infty}})$ to some $k$-additive measurable structure $(\widetilde{\widetilde{\mathcal S}}, \mathcal N_{_{\infty}})$ on $(X,G)$ having the following property:\\
there exists $\alpha _{_{0}}$$<k$, a base $\mathcal V = \{U_{_{\alpha}}\}_{_{0\leq \alpha < k}}$ of some uniformity on $\widetilde{\widetilde{\mathcal S}}$ and a strictly $k$-independent family $\{A_{i} : i\in\textsc{I}\}$ of sets from $\widetilde{\widetilde{\mathcal S}}$ with card$\hspace{.08cm}\textsc{(I)}=2^{k^{^{+}}}$ such that $(A_{i},A_{j})\notin U_{_{\alpha_{_{0}}}}$ for $i\neq j$.}
\begin{center}

\end{center}

\textbf{\underline{AUTHOR'S ADDRESS}}\

{\normalsize
\textbf{S.Basu}\\
\vspace{.02cm}
\hspace{.40cm}
\textbf{Department of Mathematics}\\
\vspace{.02cm}
\hspace{.40cm}
\textbf{Bethune College, Kolkata} \\
\vspace{.01cm}
\hspace{.40cm}
\textbf{W.B. India}\\
\vspace{.02cm}
\hspace{.1cm}
\hspace{.30cm}\textbf{{E-mail : sanjibbasu08@gmail.com}}\\
%\vspace{.0001cm}

\textbf{D.Sen}\\
\vspace{.000001cm}
\hspace{.40cm}
\textbf{Saptagram Adarsha vidyapith (High), Habra, $24$ Parganas (N)} \\
\vspace{.00001cm}
\hspace{.40cm}
\textbf{W.B. India}\\
\vspace{.00001cm}
\hspace{.1cm}
\hspace{.30cm}\textbf{{E-mail : reachtodebasish@gmail.com}}
}

\end{document}